\documentclass[final,1p,times]{elsarticle}
\usepackage{algorithm}
\usepackage{algorithmic}
\usepackage[algo2e]{algorithm2e}
\usepackage{amssymb}
\usepackage{amsthm}
\usepackage{amscd}
\usepackage{tikz}
\usepackage{amsmath}
\usepackage{amsfonts}
\usepackage{amssymb}
\usepackage{booktabs}
\usepackage{graphicx}
\usepackage{color}
\usepackage{framed}
\usepackage{comment}
\usepackage{chngcntr}
\numberwithin{equation}{section}
\newtheorem{theorem}{Theorem}[section]

\newtheorem{corollary}[theorem]{Corollary}

\usepackage{mathrsfs}
\usepackage{titletoc}
\biboptions{sort&compress}
\usepackage{subfigure}
\usepackage{hyperref}
\usepackage[pagewise,mathlines]{lineno}
\usepackage{geometry}

\usepackage{bm}

\biboptions{sort&compress}

\makeatletter

\newcommand{\Rmnum}[1]{\expandafter\@slowromancap\romannumeral #1@}
\makeatother

\usepackage{bbding}

\journal{***}

\begin{document}

\begin{frontmatter}

\title{Piecewise-linear Ricci curvature flows on weighted graphs}
		
\author[author1]{Jicheng Ma}
\ead{2019202433@ruc.edu.cn}
\author[author3]{Yunyan Yang{\footnote{corresponding author.\\ \indent\indent This research is partly supported by the National Natural Science Foundation of China (Grant No. 12271027).}}}
\ead{yunyanyang@ruc.edu.cn}

\address{$^1$School of Mathematics, Renmin University of China, Beijing, 100872, China}	

\begin{abstract}

Community detection is an important problem in graph neural networks. Recently, algorithms based on Ricci curvature flows have gained significant attention. It was suggested by Ollivier (2009), and applied to community detection by Ni et al (2019) and
Lai et al (2022). Its mathematical theory was due to Bai et al (2024) and Li-M\"unch (2025). In particular, solutions to some of these flows have existence, uniqueness and convergence. However, a unified theoretical framework has not yet been established in this field.

In the current study, we propose several unified piecewise-linear Ricci curvature flows with respect to arbitrarily selected Ricci curvatures. First, we prove that the flows have global existence and uniqueness. Second, we show that if the Ricci curvature being used is
homogeneous, then after undergoing multiple surgeries, the evolving graph has a constant Ricci curvature on each connected component.
Note that five commonly used Ricci curvatures, which were respectively defined by Ollivier, Lin-Lu-Yau, Forman, Menger and Haantjes, are all homogeneous,
and that the proof of all these results is independent of the choice of the specific Ricci curvature. Third, as an application, we apply the discrete piecewise-linear Ricci curvature flow with
surgeries to the problem of community detection. On three real-world datasets, the flow consistently outperforms baseline models and existing methods. Complementary experiments on synthetic graphs further confirm its scalability and robustness. Compared with existing algorithms, our algorithm has two advantages: it does not require curvature calculations at each iteration, and the iterative process converges.
\end{abstract}

		
\end{frontmatter}	
\section{Introduction}
In the context of Riemannian geometry, the Ricci curvature flow is a process that deforms the metric of a manifold to smooth out its curvature
\cite{Hamilton,Perelman}. Similarly, the Ricci curvature flow on weighted graphs \cite{Ollivier-1,Bai-Lin} aims to evolve the edge weights
such that the graph's structure achieves a more uniform distribution of Ricci curvature. This is based on discrete notions of Ricci curvature capturing the geometric essence of the graph. In recent years,
research on the Ricci curvature flow for weighted graphs has emerged as a field combining geometric analysis and graph theory. The adaptation of manifold-based Ricci curvature flow to the discrete setting of graphs enables the study of graphs’ geometric and topological properties through the lens of Ricci curvature evolution.

 Specifically, given a finite weighted graph $G=(V,E,\bf{w})$,
 where $V$ and $E$ represent the vertex set and edge set respectively, and
 ${\bf w}=(w_e)_{e\in E}$ denotes the edge weights, it was first observed by Ollivier  \cite{Ollivier-1}
 that the counterpart of Ricci flow on manifolds takes the form:
  \begin{equation}\label{flow}w_e^\prime(t)=-\kappa_e(t)w_e(t),\quad\forall e\in E,\end{equation}
  where $\kappa_e(t)$ is Ollivier's
  Ricci curvature for the evolving graph $G(t)=(V,E,{\bf{w}}(t))$. Shortly afterwards, Lin-Lu-Yau \cite{Lin-Lu-Yau} proposed
  a limiting version of Ollivier's Ricci curvature. Given (\ref{flow}), one might consider replacing $\kappa_e$ with Lin-Lu-Yau's Ricci curvature. Indeed, in \cite{Bai-Lin}, Bai et al. established the local existence, uniqueness, and global existence (up to surgeries) of solutions
  to (\ref{flow}) with $\kappa_e$ substituted by Lin-Lu-Yau's Ricci curvature.

   Community detection is a crucial technique in network analysis, with significant applications in sociology \cite{Scott}, biology \cite{Girvan M, Bhowmick}, and computer science \cite{Tauro S L}. Numerous algorithms have been developed for this purpose, including those presented in \cite{Clauset-Newman-Moore,Girvan M}. In 2019, Ni-Lin-Luo-Gao \cite{Ni-Lin} introduced an effective community detection method based on a discrete version of (\ref{flow}) combined with a well-designed surgical procedure. Similarly, Lai-Bai-Lin \cite{Lai X} achieved comparable results by employing a normalized Ricci flow derived from \cite{Lin-Lu-Yau, Bai-2}. Experimental results in \cite{Lai X} demonstrated the effectiveness of Lin-Lu-Yau's Ricci curvature flow for community detection. Recently, Li-M\"unch \cite{Li-Munch} proved the convergence (under surgeries) of solutions to the discrete Ollivier's Ricci flow \cite{Ni-Lin}.
   However, the convergence of solutions for the continuous Ollivier's Ricci curvature flow, the continuous Lin-Lu-Yau's Ricci curvature flow,
    or the discrete Lin-Lu-Yau's Ricci curvature flow--as applied in \cite{Lai X}--remains unresolved.

   According to \cite{Bai-Lin,Ni-Lin}, a solution to (\ref{flow}) may blow up at a finite time. However,  after applying specific
   surgical procedures,
   the solution $\textbf{w}(t)$ exists for all time $t\in[0,+\infty)$.
   Recently, in \cite{M-Y1,M-Y2}, we reformulated the Ricci curvature flow (\ref{flow}) as:
   \begin{equation}\label{rho}w_e^\prime(t)=-\kappa_e(t)\rho_e(t),\quad\forall e\in E,\end{equation}
   where $\kappa_e(t)$ denotes Lin-Lu-Yau's Ricci curvature or Ollivier's Ricci curvature on the evolving graph $G(t)=(V,E,\textbf{w}(t))$,
    and $\rho_e(t)$ represents the length (not the weight) of edge $e$. We proved that for any
   initial data $\textbf{w}(0)$, the flow
   (\ref{rho}) admits a unique global solution $\textbf{w}(t)$. Consequently, the key distinction between the two flows
   lies in the fact that global solutions to (\ref{flow}) generally require surgeries, whereas those to (\ref{rho}) do not.
   We further applied discrete versions of
   (\ref{rho}) to community detection, with our algorithms demonstrating superior performance compared to those in
   \cite{Ni-Lin,Lai X}. Nevertheless, similar to the case for
   (\ref{flow}), the convergence of solutions to (\ref{rho}) or its discrete counterpart remains an open problem.

   In the current paper, our aim is to construct a piecewise-linear Ricci curvature flow that differs from the two aforementioned flows.
   Formally, the flow is defined on each time interval $[t_{i-1},t_i)$ as
   $$w_e^{\prime}(t)=-\kappa_e(t_{i-1})w_e(t),$$
   where $\kappa_e$ denotes any established Ricci curvature, such as Ollivier's Ricci curvature,
   Lin-Lu-Yau's Ricci curvature, Forman's Ricci curvature \cite{Forman,Jost-1}, Menger's Ricci curvature \cite{Menger}, Haantjes' Ricci curvature \cite{Haantjes,Jost-2}, among others. We demonstrate that, after a finite number of surgeries, global solutions to such flows not only exist but also converge. Specifically, given an initial dataset, the piecewise-linear Ricci flow achieves constant Ricci curvature within each connected component of the evolving graph following finitely many surgeries. This result also holds for its discrete counterpart.
   The theoretical framework relies on fundamental existence theorems from classical ODE theory. Furthermore, we apply the discrete flows to community detection, with experimental results indicating that our algorithms perform as effectively as those in prior studies
   \cite{Ni-Lin,Lai X,M-Y1,M-Y2}. We also compare the performance of various Ricci curvatures in community detection tasks.

   The remainder of this paper is organized as follows. In Section \ref{Sec 2}, we propose some important definitions and main results. Section \ref{Sec 3} covers existence and convergence of solutions (Theorems \ref{P-L}), provides an overview of various Ricci curvatures, and includes an illustrative example of piecewise-linear Ricci curvature flow that incorporates these curvatures. In Section \ref{Sec 4}, we analyze  convergence of solutions under $A$-surgeries, specifically proving Theorems \ref{thm} and \ref{discrete}. Section \ref{Sec 5} outlines a
   related algorithm to community detection. Section \ref{Sec 6} presents extensive experiments evaluating the accuracy of our algorithms in addressing community detection problems. Finally, we conclude this work in Section \ref{Sec 7}.

\section{Key definitions and main results}\label{Sec 2}
Let $G=(V,E,\mathbf{w})$ be a finite weighted graph, where $V=\{z_1,z_2,\cdots,z_n\}$ is the vertex set,
$E=\{e_1,e_2,\cdots,e_m\}$ is the edge set,
$\mathbf{w}=(w_{e_1},w_{e_2},\cdots,w_{e_m})\in\mathbb{R}^m_+$ denotes the vector of edge weights,
and $\mathbb{R}^m_+$ is denoted by
$$\mathbb{R}^m_+=\{x=(x_1,x_2,\cdots,x_m)\in\mathbb{R}^m: x_i>0\, {\rm for\,\,all}\, 1\leq i\leq m\}.$$
Hereafter, we use $\kappa$ to denote any Ricci curvature on the graph $G$. Unless explicitly required, we will not distinguish between specific types of Ricci curvatures in our analysis.\\

Now, we define a continuous piecewise-linear Ricci curvature flow (\textsf{PLRF} for short) as follows. \\[1.5ex]
{\bf Definition A}. Let $0=t_0<t_1<t_2<\cdots<t_N<t_{N+1}= +\infty$ partition the time interval $[0,+\infty)$.
A function $\textbf{w}:[0,+\infty)\rightarrow\mathbb{R}^m_+$, expressed as $\textbf{w}(t)=(w_{e_1}(t),w_{e_2}(t),\cdots,w_{e_m}(t))$,
is called a continuous piecewise-linear Ricci curvature flow associated with the partition $\{t_1,t_2,\cdots,t_N\}$
if it satisfies:
\begin{itemize}
\item[$1$.]{Initial condition: $\mathbf{w}(0)=\mathbf{w}_0$;}
\item[$2$.]{ Differential equation on intervals:  For each $i=1,2,\cdots,N+1$ and each $e\in E$,
$$w_{e}^\prime(t)=-\kappa_{e}(t_{i-1})w_{e}(t)\,\,{\rm for\,\,all}\,\, t\in[t_{i-1},t_{i}),\quad
w_e|_{t_{i-1}}=w_{e}(t_{i-1}),$$
where $\kappa_e(t_{i-1})$ denotes the Ricci curvature on $e\in E$ with respect to the weighted graph
$G(t_{i-1})=(V,E,\textbf{w}(t_{i-1}))$.}
\end{itemize}

We remark that in Definition A, the finite partition $\{t_i\}_{i=1}^N$ may be replaced by an infinite partition
$\{t_k\}_{k=1}^\infty$, and that for each $i$, $\kappa_e(t_{i-1})$ may be replaced by a related real number $c_{i-1,e}$
according to specific needs.\\

Our first result is the following:

\begin{theorem}\label{P-L}
For any finite weighted graph \( G = (V, E, \mathbf{w}) \) and any partition \( \{t_1, t_2, \ldots, t_N\} \) with \( 0 =t_0< t_1 < t_2 < \cdots < t_N <t_{N+1}= +\infty \) of $[0,+\infty)$, there exists a continuous {\textsf{PLRF}} associated with \( \{t_1, t_2, \ldots, t_N\} \).
Furthermore, for each edge $e\in E$, there hold
\begin{itemize}
  \item[$1$.] if \( \kappa_e(t_N) = 0 \), then \( w_e(t) \equiv w_e(t_N) \) for all \( t \geq t_N \);
  \item[$2$.] if \( \kappa_e(t_N) \neq 0 \), then
    \[
    \lim_{t \to +\infty} w_e(t) =
    \begin{cases}
      0 & \text{if } \kappa_e(t_N) > 0, \\
      +\infty & \text{if } \kappa_e(t_N) < 0.
    \end{cases}
    \]
\end{itemize}
\end{theorem}

We conclude from Theorem \ref{P-L} that continuous {\textsf{PLRF}} always exists globally
for all $t\in [0,+\infty)$.
However, the asymptotic behavior of $\mathbf{w}(t)$ proves unsatisfactory in practical applications. To achieve refined convergence of
 $\mathbf{w}(t)$,
we define surgery as follows.\\[1.5ex]
{\bf Definition B ($A$-surgery)}.
Let $G=(V,E,\mathbf{w})$ be a finite weighted graph, and fix a real number $A>1$. For an edge $e\in E$, if
$$\frac{w_e}{\min_{e^\prime\in E^e}w_{e^\prime}}\geq A,$$
where $E^e$ denotes the set of edges in the connected component of $G$ containing $e$, then $e$ is removed from $E$.
 Let $E_A$ be the set of all such removed edges, and define the surgered graph as
$\widetilde{G}=(V,E\setminus E_A,\mathbf{w})$. The process of constructing $\widetilde{G}$ from $G$ is called an
$A$-surgery.\\

Now we propose a continuous {\textsf{PLRF}} with surgeries, which will be described in several steps.
Assume that $(t_k)_{k\in\mathbb{N}}$ is a strictly increasing sequence of positive numbers
with $t_k\rightarrow+\infty$ as $k\rightarrow+\infty$. Let
$G_0=(V,E_0,\mathbf{w}_0)$ be an initial finite weighted graph,
and $\{\kappa_{0,e}\}_{e\in E_0}$ be its Ricci curvatures on edges.
Fix a real number
\begin{equation}\label{2.1}
A > \max_{e \in E_0} \left( \frac{w_{0,e}}{\min_{e' \in E_0^e} w_{0,e'}} \right),
\end{equation}
where \( E_0^e \) denotes the connected component of \( G_0 \) containing the edge \( e \).

{\rm\textsf{Step 1 (Surgery at $t_1$).}}
Take $c_{0,e}=\kappa_{0,e}$ for each $e\in E_0$. For $t\in[0,t_1]$, the linear Ricci curvature flow
$$\left\{\begin{array}{lll}
w_e^\prime(t)=-c_{0,e}w_e(t)\\[1.2ex]
w_e(0)=w_{0,e}
\end{array}\right.$$
 admits a unique solution
$w_e(t)=w_{0,e}\exp(-c_{0,e}t)$. Denote $\mathbf{w}(t_1)=(w_e(t_1))_{e\in E_0}$. Set
$$E_{0,A}=\left\{e\in E_0: \frac{w_e(t_1)}{\min_{e^\prime\in E_0^e}w_{e^\prime}(t_1)}\geq A\right\},$$
 $E_1=E_0\setminus E_{0,A}$,
$G_1=(V,E_1, \mathbf{w}(t_1))$, and
$$c_{1,e}=\left\{\begin{array}{lll}
c_{0,e}&{\rm if}& E_1=E_0\\[1.2ex]
\kappa_e(t_1)&{\rm if}& E_1\not=E_0,
\end{array}\right.$$
where $\kappa_e(t_1)$ is the Ricci curvature on $G_1$.

{\rm\textsf{Step 2 (Surgery at $t_2$).}}
For $t\in[t_1,t_2]$, solve the linear Ricci curvature flow
$$w_e^\prime(t)=-c_{1,e}w_e(t),\quad w_e|_{t_1}=w_e(t_1),$$
yielding the unique solution
$$w_e(t)=w_{e}(t_1)\exp(-c_{1,e}(t-t_1)).$$
 Denote $\mathbf{w}(t_2)=(w_e(t_2))_{e\in E_1}$ and
$$E_{1,A}=\left\{e\in E_1: \frac{w_e(t_2)}{\min_{e^\prime\in E_1^e}w_{e^\prime}(t_2)}\geq A\right\},$$
where $E_1^e$ is the connected component of $E_1$ containing $e$. Set $E_2=E_1\setminus E_{1,A}$ and
 $G_2=(V,E_2,\mathbf{w}(t_2))$.

{\rm\textsf{Step 3}} {\rm (\textsf{Induction})}.
Suppose that we have already $G_{k-1}=(V,E_{k-1},\mathbf{w}(t_{k-1}))$ for $k\geq 2$.  Define
$$c_{k-1,e}=\left\{\begin{array}{lll}
c_{k-2,e}&{\rm if}& E_{k-1}=E_{k-2},\\[1.2ex]
\kappa_e(t_{k-1})&{\rm if}& E_{k-1}\not=E_{k-2},
\end{array}\right.$$
where $\kappa_e(t_{k-1})$ is the Ricci curvature on $G_{k-1}=(V,E_{k-1},\mathbf{w}(t_{k-1}))$.
For $t\in[t_{k-1},t_k]$, solve
$$w_e^\prime(t)=-c_{k-1,e}w_e(t),\quad w_e|_{t_{k-1}}=w_e(t_{k-1}),$$
with solution
$$w_e(t)=w_{e}(t_{k-1})\exp(-c_{k-1,e}(t-t_{k-1})).$$
Denote $\mathbf{w}(t_k)=(w_e(t_k))_{e\in E_{k-1}}$ and
$$E_{k-1,A}=\left\{e\in E_{k-1}: \frac{w_e(t_k)}{\min_{e^\prime\in E_{k-1}^e}w_{e^\prime}(t_k)}\geq A\right\},$$
where $E_{k-1}^e$ is the connected component of $G_{k-1}$ including $e$. Finally, set
$E_k=E_{k-1}\setminus E_{k-1,A}$ and $G_k=(V,E_k,\mathbf{w}(t_k))$.\\

Let us come back temporarily to discuss Ricci curvatures. We say that a Ricci curvature $\kappa$ is $\gamma$-homogeneous for
some real number $\gamma$, if for scaling weighted graphs
$G_a=(V,E,a\mathbf{w})$ with $a>0$, there hold
$$\kappa_e(G_a)=a^\gamma\kappa_e(G_1)\,\,\,{\rm for\,\,all}\,\,\, e\in E.$$
The readers will see several commonly used Ricci curvatures in Subsection \ref{Ricci_curvatures}.
It is not difficult to check that both Ollivier's Ricci curvature and Lin-Lu-Yau's Ricci curvature are $0$-homogeneous,
Forman's Ricci curvature are $1$-homogeneous, while Menger's and Haantjes' Ricci curvature are $(-1)$-homogeneous.\\

Our second result is summarized as follows:

\begin{theorem}\label{thm}
Let \( G_0 = (V, E_0, \mathbf{w}_0) \) be an initial finite weighted graph. Fix any real number $A$ satisfying (\ref{2.1}).
Let \( (t_k)_{k \in \mathbb{N}} \) be a strictly increasing sequence of positive numbers satisfying \( t_k \to +\infty \) as \( k \to +\infty \). Then, there exists a unique continuous {\textsf{PLRF}} with $A$-surgeries with respect to the partition \( (t_k)_{k\in\mathbb{N}} \).
Moreover, there exists a sufficiently large \( T > 0 \) such that:
\begin{itemize}
  \item[$1$.] No \( A \)-surgeries occur for all \( t \geq T \);
  \item[$2$.] If $\kappa$, the Ricci curvature being used, is $\gamma$-homogeneous for some real number $\gamma$, then for any $t\geq T$ and each connected component $(V^\prime,E^\prime)$ of the graph \( G(t) = (V, E(t), \mathbf{w}(t)) \), there exists a constant $\Theta=\Theta(E^\prime,\gamma,t)$ satisfying $\kappa_e(t)=\Theta$ for any $e\in E^\prime$.
\end{itemize}
\end{theorem}

However, for application, we are concerned with discrete versions of {\textsf{PLRF}} with $A$-surgeries.
This is very easy to operate: Let $(t_k)_{k\in\mathbb{N}}$ be the sequence from Theorem \ref{thm}.
A discrete {\textsf{PLRF}} with $A$-surgeries is written by
\begin{equation}\label{dis-1}\left\{\begin{array}{lll}w_e(t_k)=w_e(t_{k-1})\exp(-c_{k-1,e}(t_k-t_{k-1})),\\[1.2ex]
w_e(t_0)=w_{0,e},\,\, k=1,2,\cdots,
\end{array}\right.\end{equation}
where $c_{0,e}=\kappa_{0,e}$ and for $k\geq 1$,
\begin{equation}\label{ck}c_{k,e}=\left\{\begin{array}{lll}
c_{k-1,e}&{\rm if}& E_{k-1,A}=\varnothing\\[1.2ex]
\kappa_e(t_k)&{\rm if}& E_{k-1,A}\not=\varnothing.
\end{array}\right.\end{equation}

Our third result is stated as follows.
\begin{theorem}\label{discrete}
Let $G_0$, $A$ and $(t_k)$ be as defined in Theorem \ref{thm}.
Suppose $w_e(t_k)$ and $c_{i,e}$ are given by (\ref{dis-1}) and (\ref{ck}), respectively.
Then, there exists $\ell\in \mathbb{N}$ such that
for $k\geq \ell$, every connected
component $G^\prime(t_k)=(V^\prime,E^\prime,\mathbf{w}(t_k))$ of
the graph $G(t_k)=(V,E_{\ell-1},\mathbf{w}(t_k))$ satisfies:
\begin{itemize}
  \item [$1$.] {\bf{Edge weight ratios}:} For all $e, e^\prime\in E^\prime$,
$$\frac{w_e(t_k)}{w_{e^\prime}(t_k)}\equiv \frac{w_e(t_\ell)}{w_{e^\prime}(t_\ell)},\quad \forall k\geq \ell;$$
\item [$2$.]{\bf{Constant Ricci curvature}:} If $\kappa$, the Ricci curvature being used, is $\gamma$-homogeneous for some real number $\gamma$, then  there exists a constant $\Theta_k=\Theta_k(E^\prime,\gamma)$ such that
  $\kappa_e(t_k)=\Theta_k$ for any edge $e\in E^\prime$.
\end{itemize}
\end{theorem}

\section{Continuous \textsf{PLRF}}\label{Sec 3}

In this section, we study continuous {\textsf{PLRF}}. In particular, we first prove Theorem \ref{P-L}, where the Ricci curvature
is assumed to be arbitrary. Next, we collect several kinds of Ricci curvatures. Finally, we construct an explicit example
of continuous {\textsf{PLRF}}.

\subsection{Proof of Theorem \ref{P-L}}

 Assume $0=t_0<t_1<t_2<\cdots<t_N<t_{N+1}= +\infty$. Denote $E=\{e_1,e_2,\cdots,e_m\}$. Noticing for all $i=1,2,\cdots,N+1$,
\begin{equation}\label{linear}\left\{\begin{array}{lll}w_{e_j}^\prime(t)=-\kappa_{e_j}(t_{i-1})w_{e_j}(t)\\[1.5ex]
t\in[t_{i-1},t_i),\,\,\,j=1,2,\cdots,m\end{array}\right.\end{equation}
is an ordinary differential system with a constant coefficient matrix
  $$\mathbf{K}_{i-1}=\left(\begin{array}{ccccccc}
&\kappa_{e_1}(t_{i-1})&0&0&\cdots&0\\
&0&\kappa_{e_2}(t_{i-1})&0&\cdots&0\\
&0&0&\kappa_{e_3}(t_{i-1})&\cdots&0\\
&\vdots&\vdots&\vdots&\ddots&\vdots\\
&0&0&0&\cdots&\kappa_{e_m}(t_{i-1})
\end{array}\right),$$
we conclude from the ODE theory (\cite{Wang-Zhou-Zhu-Wang}, Chapter 6) that (\ref{linear}) has a unique solution on $[t_{i-1},t_i)$. Actually, we get the solution
\begin{equation}\label{weight}w_{e_j}(t)=w_{e_j}(t_{i-1})\exp(-\kappa_{e_j}(t_{i-1})(t-t_{i-1})),\quad\forall t\in[t_{i-1},t_i),\quad \forall j=1,2,\cdots,m.\end{equation}
Obviously $\mathbf{w}(t)=(w_{e_1}(t),w_{e_2}(t),\cdots,w_{e_m}(t))$ is continuous with respect to $t\in [0,+\infty)$. It then follows
from (\ref{weight}) that if $\kappa_{e_j}(t_N)=0$, then $w_{e_j}(t)=w_{e_j}(t_N)$ for all $t\geq t_N$; if $\kappa_{e_j}(t_N)\not=0$, then
 $$\lim_{t\rightarrow+\infty}w_{e_j}(t)=\left\{\begin{array}{lll}
0&{\rm if}& \kappa_{e_j}(t_N)>0\\[1.2ex]
+\infty&{\rm if}& \kappa_{e_j}(t_N)<0.\\[1.2ex]
\end{array}\right.$$
This completes the proof of Theorem \ref{P-L}. $\hfill\Box$

\subsection{Ricci curvatures}\label{Ricci_curvatures}

In this subsection, we will collect several Ricci curvatures on weighted finite graphs. Let
$G=(V,E,\mathbf{w})$ be a weighted finite graph.
\begin{itemize}
\item Ollivier's Ricci curvature

A function $\mu:V\rightarrow[0,+\infty)$ is said to be a probability measure if $\sum_{x\in V}\mu(x)=1$.
Let $\mu_1$ and $\mu_2$ be two probability measures. A coupling between $\mu_1$ and $\mu_2$ is defined as
a map $A:V\times V\rightarrow[0,1]$ satisfying for all $u, v\in V$,
$$\sum_{x\in V}A(u,x)=\mu_1(u),\quad\sum_{y\in V}A(y,v)=\mu_2(v).$$
The Wasserstein distance
between  $\mu_1$ and $\mu_2$ reads as
$$W(\mu_1,\mu_2)=\inf_A\sum_{u,v\in V}A(u,v)d(u,v),$$
where $A$ is taken from a set of all couplings between  $\mu_1$ and $\mu_2$.
Here and throughout,
$d(u,v)$ denotes the distance between $u$ and $v$, namely
\begin{equation}\label{distance}d(u,v)=\inf_{\gamma\in \Gamma(u,v)}\sum_{\tau\in\gamma}w_\tau,\end{equation}
$\Gamma(u,v)$ denotes the set of all paths connecting $u$ and $v$.
 Given $\alpha\in [0,1]$, an
$\alpha$-lazy one-step random walk reads as
$$\mu_{x}^\alpha(z)=\left\{\begin{array}{lll}
\alpha&{\rm if}& z=x\\[1.2ex]
(1-\alpha)\frac{w_{xz}}{\sum_{u\sim x}w_{xu}}&{\rm if}& z\sim x\\[1.2ex]
0&{\rm if}& {\rm otherwise}.
\end{array}\right.
$$
On each edge $e=xy$, Ollivier's
Ricci curvature \cite{Ollivier-1,Lin-Lu-Yau} is defined by
$$\kappa_e^\alpha=1-\frac{W(\mu_x^\alpha,\mu_y^\alpha)}{\rho_e},$$
where $\rho_e=d(x,y)$ denotes the length of $e$.
\item {Lin-Lu-Yau's Ricci curvature}

It was proved by Lin-Lu-Yau \cite{Lin-Lu-Yau} that for any fixed edge $e$, the quantity $\kappa_e^\alpha$ is concave in $\alpha\in[0,1]$
and $\kappa_e^\alpha/(1-\alpha)$ has an upper bound. As a consequence, a geometric curvature
$$\kappa_e=\lim_{\alpha\rightarrow 1}\frac{\kappa_e^\alpha}{1-\alpha}$$
is well defined.
\item {Forman's Ricci curvature}

In \cite{Forman,Jost-1}, Forman's Ricci curvature on an edge $e=xy$ was written as
$$F(e)=w_x\left(1-\sum_{e_x\sim e}\sqrt{\frac{w_e}{w_{e_x}}}\right)+w_y\left(1-\sum_{e_y\sim e}\sqrt{\frac{w_e}{w_{e_y}}}\right),$$
where $w_e$ is the weight of $e$, $w_x$ is the weight of the vertex $x$, $e_x$ denotes an edge connecting $x$,
$e_x\sim e$ means $e_x$ connects $e$ but not $e$ itself. For application, one can take
$w_x=\sum_{x\in e}w_e$ for all $x\in V$;
or
$w_x=1$ for all $x\in V$.
\item {Menger's Ricci curvature}

In \cite{Menger,Jost-2}, Menger defined a Ricci curvature on un-weighted graphs. Now we generalize Menger's Ricci curvature to
weighted graphs.
A set $T=\{x,y,z\}\subset V$ is said to be a triangle if
$y\sim x$, $y\sim z$ and $z\sim x$, i.e. $xy$, $yz$ and $zx$ are all in $E$. The distance between two vertices $u$, $v$ is written as
(\ref{distance}).
Set $a=d(x,y)$, $b=d(y,z)$, $c=d(z,x)$ and $p=(a+b+c)/2$.
Clearly we have $a\leq b+c$, $b\leq a+c$ and $c\leq a+b$, since $d(\cdot,\cdot)$ is a metric on the graph $G$.
If $T$ is a regular triangle, i.e., $a<b+c$, $b<c+a$ and $c<a+b$,
then the curvature of $T$ is defined as
$$M(T)=\frac{1}{R(T)}=\frac{\sqrt{p(p-a)(p-b)(p-c)}}{abc},$$
where $R(T)$ is the radius of the circumscribed circle of the triangle $T$. If $T$ is singular, i.e., one of the following three
alternatives holds: $a=b+c$; $b=a+c$; $c=a+b$, then the curvature of $T$ is defined as
$M(T)=0$.
Since $p\geq \max\{a,b,c\}$, the curvature of $T$ can be uniformly defined as
$$M(T)=\frac{\sqrt{p(p-a)(p-b)(p-c)}}{abc}.$$
Given an edge $e$. Let $T_e$ be a set of all triangles including $e$.
Then Menger's Ricci curvature on $e$ is defined as
$M(e)=\sum_{T\in T_e} M(T)$.
\item {Haantjes' Ricci curvature}

 Assume $\pi=x_0,x_1,\cdots,x_n$ is a simple path connecting $x_0$ and $x_n$, where each $x_{i-1}x_i$ is an edge, $i=1, 2,\cdots,n$.
The total weight of the path $\pi$ reads as
$\ell(\pi)=\sum_{i=1}^n w_{x_{i-1}x_i}$.
Haantjes' curvature \cite{Haantjes,Jost-2} on $\pi$ is defined as
$$H(\pi)=\sqrt{\frac{\ell(\pi)-d(x_0,x_n)}{d(x_0,x_n)}}\,\frac{1}{{d(x_0,x_n)}},$$
where $d(x_0,x_n)$ denotes the distance between $x_0$ and $x_n$, defined as in (\ref{distance}).
Then Haantjes' curvature on an edge $e$ is defined as
$H(e)=\sum_\pi H(\pi)$,
where $\pi$ denotes any path connecting the two vertices of $e$.
\end{itemize}

\subsection{An example of continuous {\rm{\textsf{PLRF}}}}

Take an initial graph $G_0=(V,E,\mathbf{w}_0)$, where $V=\{x_i\}_{i=1}^6$, $\mathbf{w}_0=(1,1,1,1,1,1)$, and
$E=\{x_1x_2, x_1x_3, x_2x_3,
x_2x_4, x_4x_5, x_4x_6, x_5x_6\}$. \\

\begin{figure}[H]
    \centering

\tikzset{every picture/.style={line width=0.75pt}} 

\begin{tikzpicture}[x=0.75pt,y=0.75pt,yscale=-1,xscale=1]

\draw   (79.5,149.88) -- (39.7,185.11) -- (39.3,115.12) -- cycle ;
\draw   (120.5,150.77) -- (160.03,115.24) -- (160.96,185.23) -- cycle ;
\draw    (79.5,149.88) -- (120.5,150.77) ;

\draw (29.5,100) node [anchor=north west][inner sep=0.75pt]   [align=left] {$x_1$};
\draw (29.7,185.11) node [anchor=north west][inner sep=0.75pt]   [align=left] {$x_3$};
\draw (73,151.88) node [anchor=north west][inner sep=0.75pt]   [align=left] {$x_2$};
\draw (110,151.5) node [anchor=north west][inner sep=0.75pt]   [align=left] {$x_4$};
\draw (150,100) node [anchor=north west][inner sep=0.75pt]   [align=left] {$x_5$};
\draw (150.46,185.23) node [anchor=north west][inner sep=0.75pt]   [align=left] {$x_6$};
\end{tikzpicture}
    \caption{An example of continuous \textsf{PLRF}}
    \label{fig:example PLRF}
\end{figure}
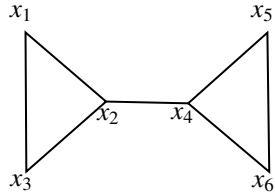

\noindent Set $t_i=0.05\times i$ for $0\leq i\leq 5$.
Let
 $$w_e(t)=\left\{\begin{array}{lll}
 w_e(t_{i-1})\exp(-\kappa_e(t_{i-1})(t-t_{i-1})),\,\, t_{i-1}\leq t<t_i\\[1.5ex]
 e\in E,\,\,1\leq i\leq 5\\[1.5ex]
 w_e(t_5)\exp(-\kappa_e(t_{5})(t-t_{5})),\,\, t\geq t_5.
 \end{array}\right.$$
Then $\mathbf{w}(t)=(w_e(t))_{e\in E}$ is the continuous \textsf{PLRF} with respect to $\{t_i\}_{i=1}^5$. It then follows that
$$w_e(t_j)=\exp\left(-\sum_{i=0}^{j-1}\kappa_e(t_i)\right),\quad 1\leq j\leq 5.$$
Following the construction and analysis of the continuous PLRF, the numerical characteristics of different curvature types across various edges in the graph can be further observed. The specific weight $w_e(t_5)$ of different curvature types for each edge are illustrated in Table \ref{PLRF example results}, where HR, MR, OR, LR, FR stand for Ricci curvatures due to Haantjes, Menger, Ollivier, Lin-Lu-Yau and
Forman respectively.
\begin{table}[htbp!]
\centering
\caption{\label{PLRF example results}Table of weight values examples for different curvature types on each edge}
\begin{tabular}{cccccccccc}
\toprule
{Curvature}\textbackslash{edge}
 &$ x_1x_2$ & $x_1x_3$ & $x_2 x_3$ & $x_2 x_4$ & $x_4x_5$ & $x_4 x_6$ & $x_5 x_6$     \\ \hline
HR                           & 0.76 & 0.76 & 0.76 & 1.00       & 0.76 & 0.76 & 0.76 \\ \hline
MR                             & 0.90 & 0.78 & 0.90 & 1.00       & 0.90 & 0.90 & 0.78 \\ \hline
OR                           & 0.91 & 0.83 & 0.91 & 1.07       & 0.91 & 0.91 & 0.83 \\ \hline
LR                         & 0.78 & 0.78 & 0.78 & 1.00       & 0.78 & 0.78 & 0.78 \\ \hline
FR                             & 3.91 & 0.78 & 3.91& \(4.32 \times 10^6\) & 3.91 & 3.91 & 0.78 \\
\bottomrule
\end{tabular}
\end{table}
We know from Table \ref{PLRF example results} that at $t_5$, the weight $w_{x_2x_4}(t_5)$ is apparently greater than
weights on other edges. As a consequence, one may delete the edge $x_2x_4$ to obtain two connected components $\{x_1,x_2,x_3\}$ and
$\{x_4,x_5,x_6\}$ of the weighted graph $G_5=(V,E\setminus\{x_2x_4\},\mathbf{w}(t_5))$. This is a simple model for community
detection.

\section{\textsf{PLRF} with $A$-surgeries}\label{Sec 4}
In this section, we concern continuous or discrete \textsf{PLRF} with $A$-surgeries, defined as in
 Section \ref{Sec 2}. Firstly we prove Theorem \ref{thm}.\\

{\it Proof of Theorem \ref{thm}.} Let
$$A>\max_{e\in G_0}\frac{w_{0,e}}{\min_{e^\prime\in G_0^e}w_{0,e^\prime}},$$
where $G_0^e$ denotes the connected component containing $e$.
Since $G_0$ is a finite graph, there are at most finitely many $A$-surgeries over time along the piecewise-linear Ricci flow.
Assume an $A$-surgery occurs at $t=t_\ell$, with no $A$-surgeries for all $t>t_\ell$. The graph after the finial $A$-surgery is denoted by $G_\ell=(V,E_\ell,\mathbf{w}(t_\ell))$. We then have
\begin{equation}\label{wei-1}w_e(t)=w_e(t_\ell) \exp\left(-\kappa_e(t_\ell)(t-t_\ell)\right),\quad\forall t\geq t_\ell,\end{equation}
 where $\kappa_e(t_\ell)$ is the Ricci curvature on $G_\ell$. For any fixed $t\in[t_\ell,\infty)$, write $\mathbf{w}(t)=(w_e(t))_{e\in E_\ell}$ and $G(t)=(V,E_\ell,\mathbf{w}(t))$.
  Let $E_\ell=\cup_{j=1}^JE^\prime_j$, where each $E^\prime_j$ is a connected component of $E_\ell$,
  and $J\geq 1$ is an integer.
  Consider any $E_j^\prime$ with at least two edges. From (\ref{wei-1}), for all $t\geq t_\ell$
  and $e, e^\prime\in E_j^\prime$,
\begin{equation}\label{ratio}\frac{w_e(t)}{w_{e^\prime}(t)}=\frac{w_e(t_\ell)}{w_{e^\prime}(t_\ell)}\exp((\kappa_{e^\prime}(t_\ell)
-\kappa_e(t_\ell))(t-t_\ell)).\end{equation}
Since no $A$-surgeries occur after $t_\ell$, we have
\begin{equation}\label{bdd}\frac{w_e(t)}{w_{e^\prime}(t)}\leq A,\,\,\forall t\geq t_\ell,\,\forall e,e^\prime\in E_j^\prime.\end{equation}
We {\it claim} that
\begin{equation}\label{claim}\kappa_e(t_\ell)=\kappa_{e^\prime}(t_\ell),\,\,\forall e,e^\prime\in E_j^\prime.\end{equation}
For otherwise, if $\kappa_e(t_\ell)<\kappa_{e^\prime}(t_\ell)$, then by (\ref{ratio}),
$w_e(t)/w_{e^\prime}(t)>2A$ for large $t$, which contradicts (\ref{bdd}); While if
$\kappa_e(t_\ell)>\kappa_{e^\prime}(t_\ell)$, we have by (\ref{ratio}) that for sufficiently large $t$,
$$\frac{w_{e^\prime}(t)}{w_{e}(t)}=\frac{w_{e^\prime}(t_\ell)}{w_{e}(t_\ell)}\exp((\kappa_{e}(t_\ell)
-\kappa_{e^\prime}(t_\ell))(t-t_\ell))>2A,$$
which also contradicts (\ref{bdd}). Thus the claim follows.

Fix any $j$, $1\leq j\leq J$. Clearly, substituting (\ref{claim}) into (\ref{ratio}) gives
$$
\frac{w_e(t)}{w_{e^\prime}(t)}=\frac{w_e(t_\ell)}{w_{e^\prime}(t_\ell)},\,\,\,\forall t\geq t_\ell,\,\forall e,e^\prime\in E_j^\prime.
$$
Since the Ricci curvature being used is $\gamma$-homogeneous, we conclude from (\ref{wei-1}) that
$$\kappa_e(t)=\kappa_e(t_\ell)\exp{(-\gamma \kappa_e(t_\ell)(t-t_\ell))},\quad\forall t\geq t_\ell.$$
This together with (\ref{claim}) implies
 $$\Theta=\kappa_e(t_\ell)\exp{(-\gamma \kappa_e(t_\ell)(t-t_\ell))}$$
 is a constant depending only on $E_j^\prime$, $\gamma$ and $t$, and thus completes the proof of the theorem.
 $\hfill\Box$\\

Secondly we prove Theorem \ref{discrete}.\\

{\it Proof of Theorem \ref{discrete}.} Since $G_0$ is a finite graph, there exists some $\ell\in \mathbb{N}$ such that $E_\ell\not=E_{\ell-1}$ and $E_{k}=E_\ell$ for all $k\geq \ell$. This is equivalent to saying $t_\ell$ is the time of the last surgery.  Hence for all $k\geq \ell$, there hold $c_{k,e}=\kappa_e(t_\ell)$ and
$$w_e(t_k)=w_e(t_\ell) \exp(-\kappa_e(t_\ell)(t_k-t_\ell)).$$
Since the remaining part is completely analogous to the proof of Theorem \ref{thm}, we leave the details to the interested
readers. $\hfill\Box$\\

As a consequence of Theorems \ref{thm} and \ref{discrete}, we have the following:

\begin{corollary}
Let $G_0=(V,E_0,\mathbf{w}_0)$, $A$, $(t_k)_{k\in\mathbb{N}}$ and $\mathbf{w}(t)$ be as in Theorems \ref{thm} and \ref{discrete} respectively.
If $\kappa$, the Ricci curvature being used, is Ollivier's Ricci curvature or Lin-Lu-Yau's Ricci curvature, then there exists some $t_\ell>0$ such that for all $t\geq t_\ell$ ($t_k\geq t_\ell$), each connected component $(V^\prime,E^\prime)$ of the graph $G(t)=(V,E(t),\mathbf{w}(t))$ has a uniform constant Ricci curvature. In particular, for all
$t\geq t_\ell$ ($t_k\geq t_\ell$), there hold
$$\kappa_e(t)=\kappa_e(t_\ell)=\kappa_{e^\prime}(t_\ell)\,\,(\kappa_e(t_k)=\kappa_e(t_\ell)=\kappa_{e^\prime}(t_\ell))\,\,
{\rm for\,\,all}\,\, t\geq t_\ell\,\, (t_k\geq t_\ell)\,\, {\rm and\,\,all}
\,\, e,e^\prime\in E^\prime.$$
\end{corollary}

{\it Proof.} We only prove the case of continuous {\textsf{PLRF}} with $A$-surgeries, since the discrete case is almost the same.
Let $\kappa$ be Ollivier's Ricci curvature or Lin-Lu-Yau's Ricci curvature (Subsection \ref{Ricci_curvatures}). Clearly, $\kappa$ is scaling-invariant: for any graph
$\tilde{G}=(\tilde{V},\tilde{E},\tilde{\mathbf{w}})$ and the scaling graph $\tilde{G}_a=(\tilde{V},\tilde{E},a\tilde{\mathbf{w}})$
with some constant $a>0$, one has
\begin{equation}\label{scaling}\kappa_e(\tilde{G}_a)=\kappa_e(\tilde{G}),\quad\forall e\in E.\end{equation}
 As in the proof of Theorem \ref{thm}, one finds the time $t_\ell$ of the final $A$-surgery. Moreover
 $$w_e(t)=w_e(t_\ell)\exp{(-\kappa_e(t_\ell)(t-t_\ell))},\quad\forall t\geq t_\ell,\,\,\forall e\in E_\ell.$$
 For each connected component $E^\prime$ of $E_\ell$, it follows from (\ref{claim}) that
 $$\kappa_e(t_\ell)=\kappa_{e^\prime}(t_\ell),\quad\forall e,e^\prime\in E^\prime.$$
 Since there is no surgery after $t_\ell$, one understands that $E^\prime$ is also a connected component of $E(t)$ with $t\geq t_\ell$.
 Setting $\kappa_e(t_\ell)\equiv c$ for all $e\in E^\prime$, one gets
 $$w_e(t)=w_e(t_\ell)\exp(-c(t-t_\ell)),\quad\forall t\geq t_\ell,\,\,\forall e\in E^\prime.$$
 This together with (\ref{scaling}) leads to
 $$\kappa_e(t)=\kappa_e(t_\ell)=c,\quad \forall t\geq t_\ell,\,\,\forall e\in E^\prime,$$
 and thus completes the proof of the corollary. $\hfill\Box$\\

In the following three sections, we still denote the discrete {\textsf{PLRF}} with $A$-surgeries as {\textsf{PLRF}} for simplicity.

\section{Algorithm design for applying {\textsf{PLRF}} to community detection}\label{Sec 5}

Just as Ricci curvature flow based on Ollivier's Ricci curvature or Lin-Lu-Yau's Ricci curvature \cite{Ni-Lin,Lai X,M-Y1,M-Y2},
the {\textsf{PLRF}} can also be applied to the community detection problem. In particular, we utilize Theorem \ref{discrete} to design a
pseudo code of our {\textsf{PLRF}} algorithm, recorded as in Algorithm 1.

\begin{algorithm}[H]
\SetAlgoLined
\KwIn{ an undirected finite network \( G = \left( {V,E_0,w_0}\right)\), threshold \(A\), time series \(0=t_0 <t_1<t_2< \ldots< t_N\) }
\KwOut{community detection  results of $G$}
\(i \gets 0\)\\
\While{\(i<N\)}{
    \If{\(i = 0\)}{
        \(t_i \gets t_1\)\\
        \For{\(e \in E_0\)}{
            \(c_{0,e} \gets \kappa_{0,e}\) \\
            \(w_e(t) \gets w_{0,e} \cdot \exp(-c_{0,e} \cdot t)\)
        }
        \(\mathbf{w}(t_i) \gets (w_e(t_i))_{e\in E_0}\)\\
        \(E_{i,A} \gets \{e\in E_0: \frac{w_e(t_i)}{\min_{e'\in E_0^e} w_{e'}(t_i)} \geq A\}\)\\
        \(E_{i + 1} \gets E_0 \setminus E_{i,A}\)\\
        \(G_{i + 1} \gets (V, E_{i + 1}, \mathbf{w}(t_i))\)
    }
    \Else{
        \For{\(e \in E_i\)}{
            \If{\(E_i = E_{i - 1}\)}{
                \(c_{i,e} \gets c_{i - 1,e}\)
            }
            \Else{
                \(c_{i,e} \gets \kappa_e(t_{i - 1})\)
            }
        }
        \(t \gets t_{i+1}\)\\
        \For{\(e \in E_i\)}{
            \(w_e(t) \gets w_e(t_i) \cdot \exp(-c_{i,e} \cdot (t - t_i))\)
        }
        \(\mathbf{w}(t_{i + 1}) \gets (w_e(t_{i + 1}))_{e\in E_i}\)\\
        \(E_{i,A} \gets \{e\in E_i: \frac{w_e(t_{i + 1})}{\min_{e'\in E_i^e} w_{e'}(t_{i + 1})} \geq A\}\)\\
        \(E_{i + 1} \gets E_i \setminus E_{i,A}\)\\
        \(G_{i + 1} \gets (V, E_{i + 1}, \mathbf{w}(t_{i + 1}))\)
    }

    \(i \gets i + 1\)
}
compute connected components $C_1\cup \dots \cup C_k$ of $G_N$\\
\For{\(i\gets1\) \textbf{to} \(|V|\)}{\If{\(v_i\in C_j\)}{set clustering labels \(Y_i=j\)}}
calculate the accuracy of community detection

\Return the accuracy of community detection
\caption{Community detection using {\textsf{PLRF}}}
\end{algorithm}

This algorithm features an outer \(\texttt{while}\) loop that executes \(N\) times. For each iteration, the dominant computational cost arises from the curvature calculation, which has a complexity of \(O(|E|D^3)\), where $|E|$ is the number of edges and \(D\) is the average degree. When \(i = 0\), a \(\texttt{for}\) loop iterates over edges in \(E_0\), with the step \(c_{0,e} \gets \kappa_{0,e}\) contributing the \(O(|E|D^3)\) term, while other operations (e.g., weight updates and set operations) are linear in the number of edges but negligible compared to the curvature computation. For \(i > 0\), nested loops over \(E_i\) again involve curvature calculations \(c_{i,e} \gets \kappa_e(t_{i-1})\), each incurring \(O(|E|D^3)\). Since the outer loop runs \(N\) times and each iteration is dominated by the \(O(|E|D^3)\) curvature step, the total complexity of the main loop is \(O(N|E|D^3)\). Post-processing steps (computing connected components and assigning cluster labels) have lower-order complexities \(O(|V| + |E_N|)\) and \(O(|V|)\), respectively, which are secondary compared to the dominant term. Thus, the overall time complexity of the algorithm is governed by \(O(N|E|D^3)\).

\section{Experiments}\label{Sec 6}

In this section, we present the datasets, baseline methods, and metrics used in our experiments. We validate the discrete {\textsf{PLRF}} algorithm (still denoted by \textsf{PLRF} for short) by comparing its performance to other methods on real-world networks--such as Zachary's Karate Club, College Football, and Facebook datasets--as well as synthetic LFR benchmarks of varying sizes. Noise levels are systematically varied during these comparisons to assess robustness. The code is available at https://github.com/mjc191812/Piecewise-linear-Ricci-curvature-flows-on-weighted-graphs.

\subsection{Real datasets and synthetic datasets}
 For the real-world datasets, we select three distinct scale community graphs to evaluate the performance of the {\textsf{PLRF}} on real networks. Basic information for real-world networks is listed in Table \ref{tab:1}. The Zachary's Karate Club dataset \cite{Zachary} is a classic social network analysis benchmark consisting of 34 nodes (representing club members) and 78 undirected edges (representing interactions). The ground-truth community structure of this network is well-documented, reflecting two distinct factions that emerged due to leadership disputes. The College Football dataset \cite{Girvan M} models the 2000 NCAA Division football season, containing 115 nodes and 613 undirected edges. The vertices correspond to the teams, while the edges represent the matches between the teams. The Facebook Network dataset is a real-world social graph crawled from the Stanford Network Analysis Project (SNAP)\cite{Jure L}. Its benchmark community structure is defined by explicit attributes such as academic departments, interests, and social affiliations, making it ideal for evaluating community detection algorithms in complex social systems.
\begin{table}[H]
\centering
\caption{\label{tab:1}Summary of real-world network characteristics.}
\begin{tabular}{ccccccc}
\toprule
networks & vertexes & edges & \#Class & AvgDeg & density  &Diameter\\
\midrule
Karate   & 34      & 78    & 2                     & 4.59   & 0.139  &5 \\
Football & 115     & 613   & 12                    & 10.66  & 0.094   &4\\
Facebook & 775     & 14006 & 18                    & 36.15  & 0.047    &	8\\
\bottomrule
\end{tabular}
\end{table}

For synthetic datasets, we employed LFR benchmark networks \cite{Lancichinetti} that feature well-defined community structures. Table \ref{tab:lfr_parameters} outlines the key parameters used during network generation, where \(\mu\) represents the inter-community connection probability (ranging between 0 and 1). Higher \(\mu\) values indicate weaker community structures. These synthetic datasets serve as a controlled testing environment, enabling systematic evaluation of algorithm performance across varying noise levels. Three network series with distinct scales were generated using the parameters listed in Table \ref{tab:lfr_parameters_setting}. Testing scenarios incorporated \(\mu\) values ranging from 0.1 to 0.8, with ten networks generated per \(\mu\) value to mitigate random variability. Final results report mean performance metrics, providing a comprehensive demonstration of the algorithm’s stability and noise resilience.

\begin{table}[htbp]
  \centering
  \caption{Main parameters of LFR benchmark network}
  \begin{tabular}{cc}
    \toprule
    Parameter & Meaning \\
    \midrule
    $|V|$ & Number of nodes in the network \\
    ave\_degree & Average degree of the network \\
    max\_degree & Maximum degree of the network \\
    min\_C & Minimum number of nodes in a community \\
    max\_C & Maximum number of nodes in a community \\
    \(\mu\) & Probability of a node connecting to the outside of the community \\
    \bottomrule
  \end{tabular}
  \label{tab:lfr_parameters}
\end{table}

\begin{table}[htbp]
  \centering
  \caption{The parameter settings of the LFR benchmark generator.}
    \begin{tabular}{lcccccc}
    \toprule
    Network & $|V|$ &  ave\_degree & max\_degree &  min\_C & max\_C & $\mu$ \\
    \midrule
    LFR500  & 500   & 20             & 50         & 10           & 50           & 0.1--0.8 \\
    LFR1000 & 1000  & 20             & 50         & 10           & 50           & 0.1--0.8 \\
    LFR5000 & 5000  & 20             & 50         & 10           & 50           & 0.1--0.8 \\
    \bottomrule
    \end{tabular}
     \label{tab:lfr_parameters_setting}
\end{table}

\subsection{Evaluation and comparison algorithms}
 We will use three metrics to assess community detection's precision in real-world datasets. The normalized mutual information (NMI) \cite{Danon-Guilera-Duch} are chosen as the criteria for evaluating the quality of clustering accuracy when compared to the ground truth. Furthermore, modularity (Q) \cite{Clauset-Newman-Moore, M.Newman} is chosen to measure the robustness of the community structure of a given graph without relying on ground-truth clustering.
To be more specific,  we let
\(\{ U_1, U_2, \ldots, U_I \} \) and \(\{ W_1, W_2, \ldots, W_J \} \) be two partitions of the set \( S \) of $n$ vertices (nodes). Denote \( m_{ij} =|U_i\cap W_j|\) the number of vertices in \( U_i\cap W_j \), while \( c_i \) and \( d_j \) represent the numbers of vertices in \( U_i \) and \( W_j \), respectively. All these quantities are listed in Table \ref{tab:2}.
\begin{table}[H]
\centering
\caption{\label{tab:2}Contingency table for community detection metrics.}
\begin{tabular}{c|cccc|c}
$U \backslash W$         & $W_1$    & $W_2$    & $\cdots$ & $W_J$    & \text{sums} \\ \hline
$U_1$                    & $m_{11}$ & $m_{12}$ & $\cdots$ & $m_{1J}$ & $c_1$                    \\
$U_2$                    & $m_{21}$ & $m_{22}$ & $\cdots$ & $m_{2J}$ & $c_2$                    \\
$\vdots$                 & $\vdots$ & $\vdots$ & $\ddots$ & $\vdots$ & $\vdots$                 \\
$U_I$                    & $m_{I1}$ & $m_{I2}$ & $\cdots$ & $m_{IJ}$ & $c_I$                    \\ \hline
\text{sums} & $d_1$    & $d_2$    & $\cdots$ & $d_J$    &
\end{tabular}
\end{table}

 Then the explicit expressions of the above mentioned two criteria are written below.
\begin{itemize}
    \item \textbf{Normalized mutual information}
    \[
        \text{NMI} = \frac{-2 \sum_{i=1}^I \sum_{j=1}^J m_{ij} \log \left( \frac{m_{ij} n}{c_i d_j} \right)}{\sum_{i=1}^I c_i \log \left( \frac{c_i}{n} \right) + \sum_{j=1}^J d_j \log \left( \frac{d_j}{n} \right)}.
        \]
   NMI measures the similarity between two clusterings by computing the mutual information between the two clusterings and the mutual information between the clustering and the ground-truth partition. The higher the NMI score, the better the clustering results.
    \item \textbf{Modularity}
    \[
    Q = \sum_{k=1}^M \left( \frac{C_k}{|E|} - \beta \left( \frac{D_k}{2|E|} \right)^2 \right),
    \]
     where \( M \) represents the number of communities, \( C_k \) is the number of connections within the $k$th community, \( D_k \) is the total degree of vertices in the $k$th community, and \( \gamma \) is a resolution parameter, with a default value of \(1\). The value of \( Q \) ranges from \(-0.5\) to \(1\). Modularity measures the quality of a partition of a graph based on the
degree of its connectivity. The higher the modularity score, the better the partition
results.
\end{itemize}

We used all the aforementioned networks as inputs for our experiments, comparing our method against several classical and deep learning approaches. The algorithms selected for comparison include: Girvan-Newman \cite{Girvan M}, GraphSage \cite{Xu-Ruan}, Infomap \cite{Rosvall}, Louvain \cite{Blondel}, Label Propagation Algorithm (LPA) \cite{Raghavan-Albert}, VGAE \cite{Kipf-Welling}, and Walktrap \cite{Pons}. We applied the
{\textsf{PLRF}} method with a hyperparameter setting of
$$A=2\max_{e\in G_0}\frac{w_{0,e}}{\min_{e^\prime\in G_0^e}w_{0,e^\prime}}.$$

\subsection{Results and analysis}
\subsubsection{The results for real-world data}
To validate the effectiveness of {\textsf{PLRF}} in community detection, experiments were conducted on both real-world and synthetic datasets, followed by comparisons with some popular and advanced algorithms. Table \ref{tab:results on real datasets} presents the NMI and Q values of {\textsf{PLRF}} (using Ollivier's Ricci curvature) and other algorithms on real-world datasets. The largest values of the two indexes on each network are typed in bold.

\begin{figure}[htbp]
    \centering
    \subfigure[The ground-truth community structure]{%
        \includegraphics[width=0.48\textwidth]{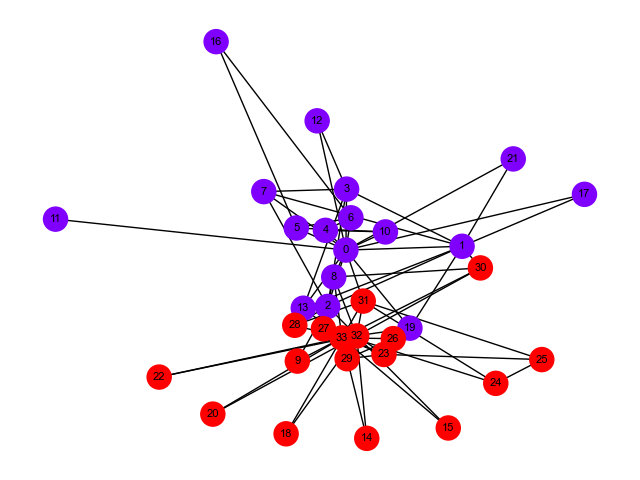}%
    }
    \hfill
    \subfigure[The community structure detected by PLRF]{%
        \includegraphics[width=0.48\textwidth]{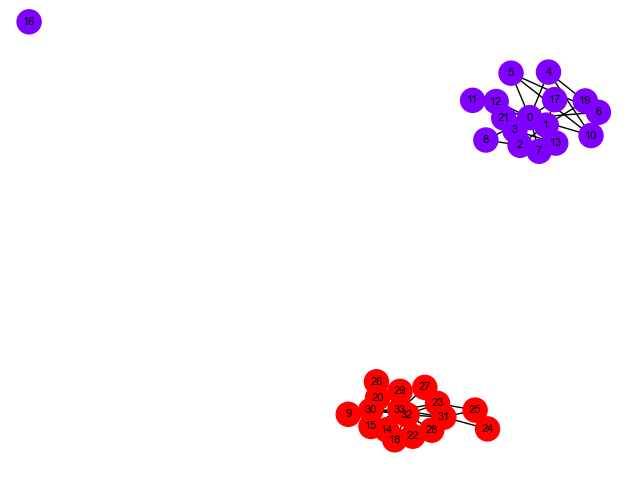}%
    }
    \caption{Community detection on the Karate club network of PLRF.}
    \label{fig:karate_example}
\end{figure}

\begin{table}[htbp]
  \centering
  \caption{NMI and Modularity on real datasets.}
  \label{tab:results on real datasets}
  \begin{tabular}{l cc cc cc}
    \toprule
    Network       & \multicolumn{2}{c}{Karate} & \multicolumn{2}{c}{Football} & \multicolumn{2}{c}{Facebook} \\
    \cmidrule(lr){2-3} \cmidrule(lr){4-5} \cmidrule(lr){6-7}
    Methods              & NMI     & Q      & NMI     & Q      & NMI          & Q      \\
    \midrule
    Girvan Newman & 0.73    & 0.48   & 0.36    & 0.50   & 0.16         & 0.01   \\
    GraphSage     & 0.74    & 0.38   & 0.30    & 0.53   & 0.30         & 0.29   \\
    Infomap       & 0.51    & 0.44   & 0.58    & 0.01   & \textbf{0.75} & 0.30   \\
    Louvain       & 0.38    & 0.39   & 0.48    & 0.55   & 0.52         & 0.45   \\
    LPA           & 0.36    & 0.54   & 0.87    & 0.90   & 0.65         & 0.51   \\
    {\textsf{PLRF}}          & \textbf{0.93} & \textbf{0.61} & \textbf{0.94} & \textbf{0.92} & 0.72         & \textbf{0.95} \\
    VGAE          & 0.61    & 0.51   & 0.69    & 0.55   & 0.51         & 0.44   \\
    Walktrap      & 0.49    & 0.01   & 0.88    & 0.01   & 0.72         & 0.30   \\
    \bottomrule
  \end{tabular}
\end{table}
The comparative analysis of community detection methodologies across three representative network datasets (Karate, Football, and Facebook) highlights the distinguishing characteristics of the proposed {\textsf{PLRF}} approach. On the Karate and Football networks, {\textsf{PLRF}} demonstrates exceptional performance in detecting ground-truth community structures, achieving state-of-the-art NMI scores of 0.93 and 0.94 respectively. These results represent significant improvements of 25.6\% and 6.8\% over the second-best baseline methods, underscoring its effectiveness in moderately sized networks. Figure \ref{fig:karate_example} shows the community division results of
{\textsf{PLRF}} on the Karate network. After the iteration and surgery, the network is divided into two communities (red and blue), which is completely consistent with the real structure.

Notably, {\textsf{PLRF}} maintains consistent modularity optimization across all evaluated datasets. It achieves the highest Q scores: 0.61 for Karate, 0.92 for Football, and 0.95 for Facebook. The particularly strong performance on the Facebook network is remarkable, where it attains near-perfect modularity while maintaining competitive NMI (0.72). This balance suggests an ability to preserve both structural cohesion and functional separation in large-scale networks.

The method demonstrates robustness across varying network scales. While NMI scores naturally decrease with increasing network size (from 0.93 in small-scale Karate to 0.72 in large-scale Facebook), modularity remains consistently high. This pattern aligns with known scaling challenges in community detection, where information-theoretic metrics are more sensitive to network size than structural measures.

Comparative evaluations against alternative approaches reveal distinct advantages. {\textsf{PLRF}} outperforms deep learning methods (GraphSage, VGAE) by 25.6-52.4\% in NMI, indicating limitations of neural methods in preserving community structures. Compared to optimization-based methods (Louvain, Infomap), it achieves 38.2-94.7\% higher modularity scores, highlighting superior capability in modularity maximization.

Anomalous results on the Facebook network require further investigation. While {\textsf{PLRF}} achieves high modularity (0.95), its NMI (0.72) lags slightly behind Infomap's 0.75. Potential contributing factors include variations in network density (Facebook: 0.047 vs. Football: 0.094), differences in ground-truth community granularity, and edge sparsity patterns affecting information-theoretic metrics.

\begin{figure}[htbp]
    \centering
    \includegraphics[width=1\textwidth]{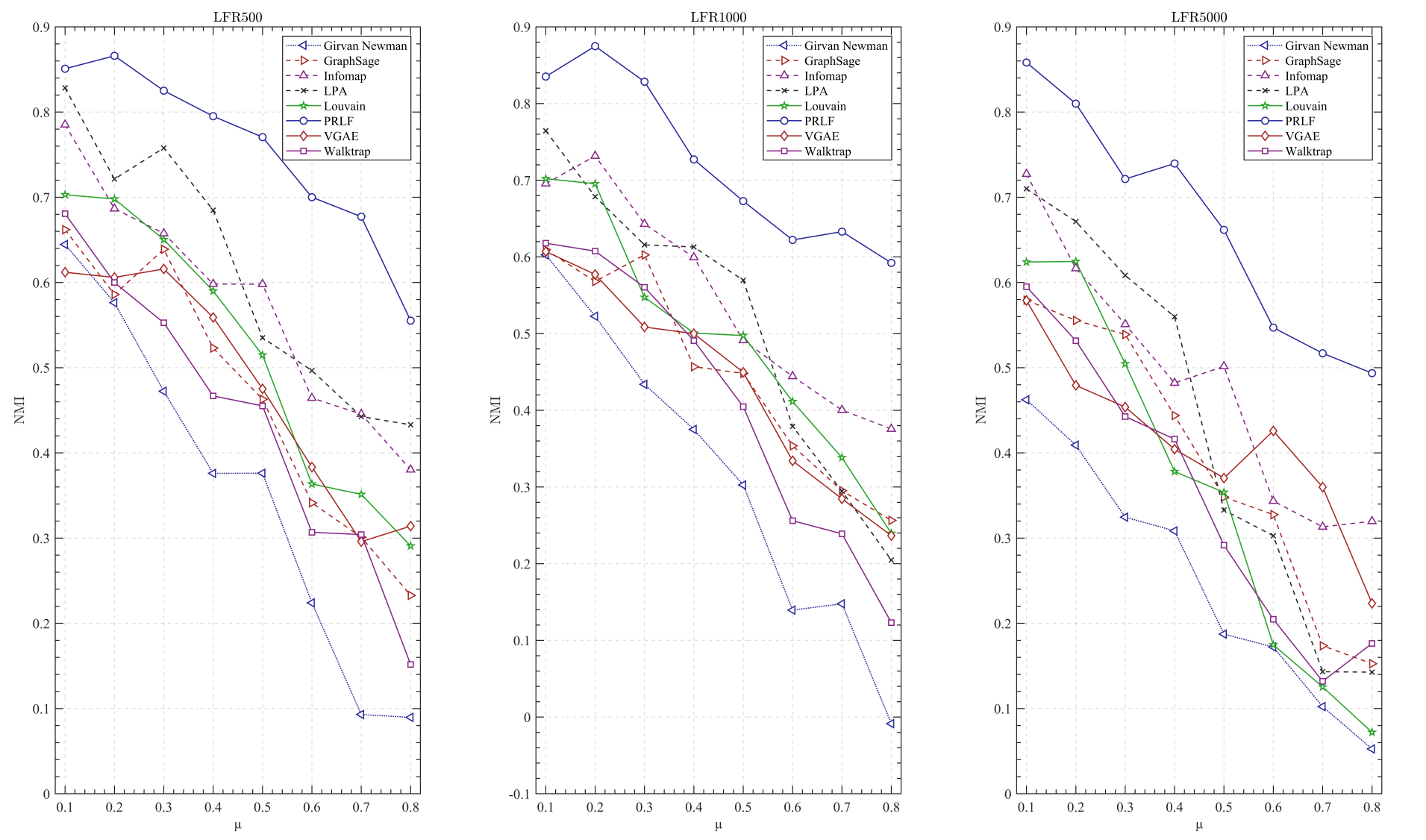}
    \caption{The NMI on the artificial networks.}
    \label{fig:NMI on synthetic network}
\end{figure}

\begin{figure}[htbp]
    \centering
    \includegraphics[width=1\textwidth]{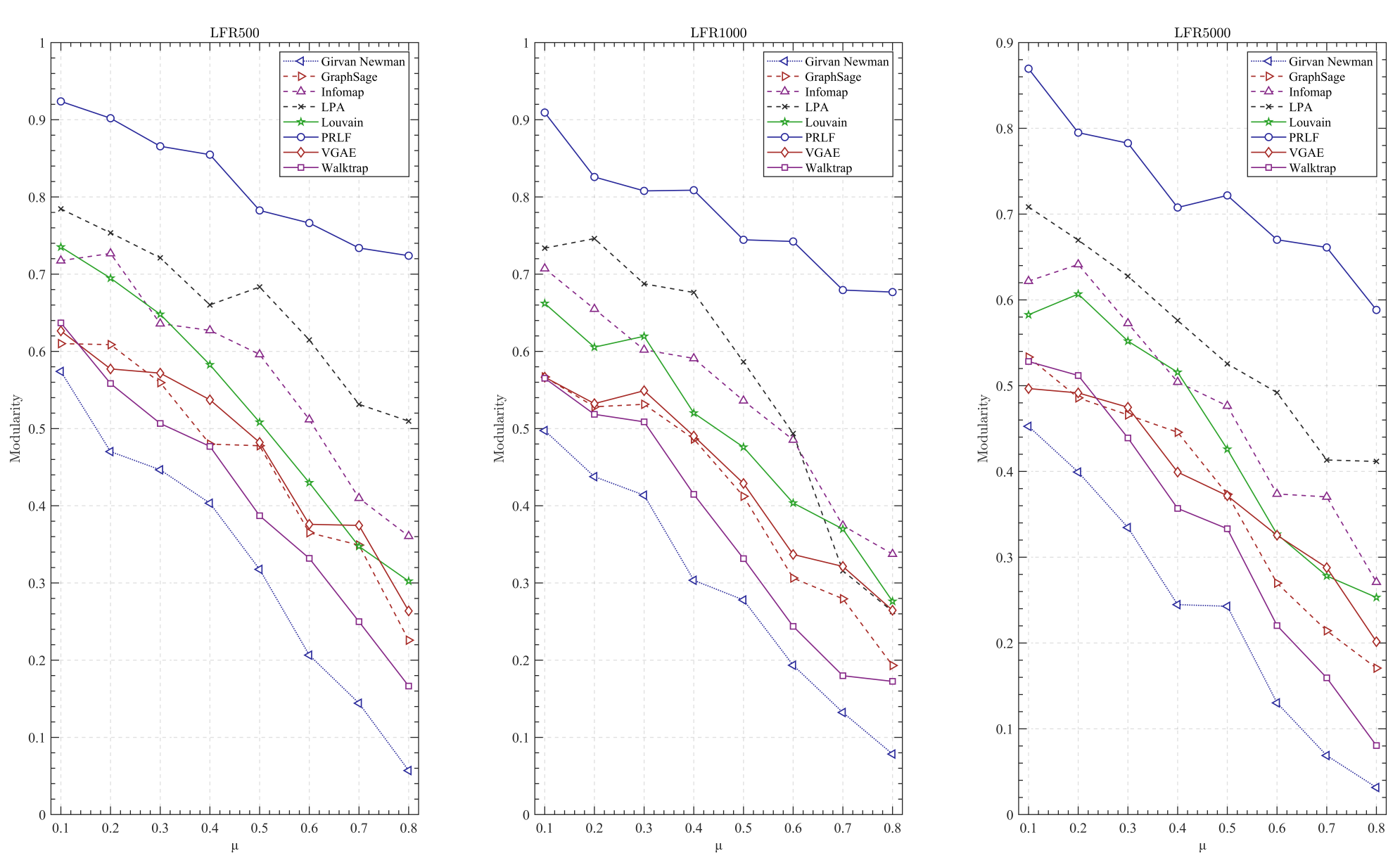}
    \caption{The Modularity on the artificial networks.}
    \label{fig:Q on synthetic network}
\end{figure}

\subsubsection{The results for synthetic data}
Figures \ref{fig:NMI on synthetic network} and \ref{fig:Q on synthetic network} demonstrate the performance comparison on LFR synthetic datasets.
For synthetic datasets, the proposed method {\textsf{PLRF}} achieves the best performance on all the evaluated networks, and its superiority extends to $\mu \leq$ 0.8 in the smaller-scale series (e.g., LFR500), $\mu \leq$ 0.8 in the medium-scale series (LFR1000), and $\mu \leq$ 0.8 in the larger-scale series (LFR5000). On the networks within each series, {\textsf{PLRF}}’s rank in terms of NMI consistently remains among the top one. Its performance is less satisfactory only on a very big subset of networks with $\mu = 0.8$, where the network structure approaches that of a random graph, rendering the community boundaries ambiguous. Even in these cases, however, the NMIs detected by {\textsf{PLRF}} are still significantly higher than those obtained by the competing algorithms.

As the difficulty level of community detection increases ($\mu$ values increase), the performance of all algorithms declines. However,
{\textsf{PLRF}} algorithm shows a relatively smaller decline in performance, highlighting its strong robustness. For example, in the moduality of the LFR5000 dataset, as $\mu$ increases from 0.1 to 0.8, while algorithms like Girvan Newman and VGAE experience a sharp drop in modularity, the decline of {\textsf{PLRF}} is more moderate. Starting from a high modularity value close to 0.9 at $\mu$ = 0.1, it still maintains a relatively high value compared to many other algorithms even when $\mu$ reaches 0.8. Similarly, in the NMI of the LFR500 dataset, as $\mu$ increases, other algorithms show a more significant decrease in NMI values. In contrast, {\textsf{PLRF}}'s NMI value drops more gradually, which shows that
{\textsf{PLRF}} can maintain a relatively stable performance in the face of increasing noise and complexity in the network structure. On the Facebook network, {\textsf{PLRF}} achieves a high modularity (Q=0.95) but slightly lower NMI (0.72) compared to Infomap (0.75). This discrepancy arises because modularity rewards dense internal connections, whereas NMI penalizes misclassifications in large-scale networks. {\textsf{PLRF}}’s curvature-based flow prioritizes structural coherence over strict alignment with ground-truth labels, leading to a trade-off between topological consistency and information-theoretic accuracy.

A thorough analysis of the experimental results on both real and virtual datasets reveals that no single algorithm consistently outperforms the others across all types of networks. However, the proposed {\textsf{PLRF}} method consistently demonstrates the ability to extract high-quality community structures, especially in networks with moderate mixing parameters. The experimental results confirm that {\textsf{PLRF}} excels in both NMI and the modularity of community detection.

The performance superiority of {\textsf{PLRF}} can be attributed to its innovative integration of discrete Ricci curvature principles, which effectively capture essential topological structural information in complex networks. The results demonstrate that Ricci curvature better characterizes network functional hierarchies than pure connectivity patterns. Our implementation extends these insights through adaptive curvature thresholding, enabling automatic detection of scale-dependent community structures. While ensuring theoretical convergence, piecewise-linear iteration also greatly improves the actual calculation speed.

\subsubsection{The ablation study for the effects of piecewise-linear design}
 We conduct an ablation study to demonstrate the effects of piecewise linear. We compared {\textsf{PLRF}} with other algorithms
  \cite{Lai X,M-Y1,M-Y2,Ni-Lin} based on the Ricci curvature flow. The results are summarized in Table \ref{tab:ablation_piecewise_linear}.
\begin{table}[htbp]
  \centering
\caption{Ablation study of piecewise-linear Ricci flow-based community detection}
\label{tab:ablation_piecewise_linear}
  \begin{tabular}{l cc cc cc}
    \toprule
    Network       & \multicolumn{2}{c}{Karate} & \multicolumn{2}{c}{Football} & \multicolumn{2}{c}{Facebook} \\
    \cmidrule(lr){2-3} \cmidrule(lr){4-5} \cmidrule(lr){6-7}
    Methods              & NMI     & Q      & NMI     & Q      & NMI          & Q      \\
    \midrule
    PLRF          & \textbf{0.93} & 0.61 & \textbf{0.94} & \textbf{0.92} & 0.72         & \textbf{0.95} \\
DORF                                             & 0.57 & 0.69          & \textbf{0.94} & 0.91          & \textbf{0.73} & 0.68          \\
NDORF                                            & 0.57 & 0.69          & \textbf{0.94} & 0.91          & \textbf{0.73} & 0.68          \\
NDSRF                                            & 0.57 & 0.68          & \textbf{0.94} & 0.91          & \textbf{0.73} & 0.68          \\
Rho                                              & 0.68 & 0.82          & 0.92          & 0.90          & 0.72          & 0.63          \\
RhoN                                             & 0.68 & \textbf{0.84}          & 0.93          & \textbf{0.92} & 0.72          & \textbf{0.95}          \\
    \bottomrule
  \end{tabular}
\end{table}

{\textsf{PLRF}} consistently achieves the highest or near-highest  values of NMI and Q across three real-world networks. Although RhoN achieves a slightly higher Q (0.84) than {\textsf{PLRF}} (0.61), the very large gain in NMI suggests that {\textsf{PLRF}} trades a modest drop in modularity for a substantial improvement in label agreement. {\textsf{PLRF}} obtains the top NMI scores on the Karate (0.93) and Football (0.94) networks and ties or outperforms all baselines in modularity, demonstrating that the piecewise-linear adjustment substantially enhances both accuracy and community quality.

Introducing a piecewise-linear component into the Ricci curvature flow framework provides a significant advantage. {\textsf{PLRF}} demonstrates superior or competitive modularity in all cases, and it leads in NMI in two out of three networks. This indicates that piecewise linear scaling of curvature better aligns the flow dynamics with the community structure, especially in small to medium-sized graphs, while still performing well on large-scale social networks. The piecewise-linear design also offers faster calculation speed and improved convergence guarantees.

\subsubsection{The ablation study for effects of curvature type}

We conducted an ablation study to examine the impact of different types of curvature. Specifically, we compared {\textsf{PLRF}} using Ollivier's Ricci curvature with other Ricci curvatures. The results, summarized in Table \ref{tab:ablation_curvature_type}, demonstrate that the choice of Ricci curvature significantly affects both detection quality and computational cost (measured in runtime in seconds). The best indicators are highlighted in bold, and the shortest run-time is underlined. OOT (out-of-time) thresholds are defined as failure on this GPU or exceeding 24 hours on an Intel i9-12900KF CPU with 16 cores.
\begin{table}[htbp]
  \centering
\caption{Ablation study on the impact of curvature type in {\textsf{PLRF}} community detection}
\label{tab:ablation_curvature_type}
  \begin{tabular}{l cc c cc c cc c}
    \toprule
    Network       & \multicolumn{3}{c}{Karate}  & \multicolumn{3}{c}{Football}  & \multicolumn{3}{c}{Facebook}  \\
    \cmidrule(lr){2-4}  \cmidrule(lr){5-7}  \cmidrule(lr){8-10}
    Curvatures              & NMI     & Q      & Time & NMI     & Q      & Time & NMI          & Q      & Time \\
    \midrule
    Ollivier           & \textbf{0.93} & 0.61 & 2.55 & \textbf{0.94} & \textbf{0.92} & 8.06 & 0.72         & \textbf{0.95} & 1775.81 \\
Lin-Lu-Yau                                              & 0.57 & \textbf{0.86} & 0.63 & 0.93 & 0.91 & 5.61 & 0.71 & 0.93 & 1396.94 \\
Forman                                             & 0.49 & 0.67 & 0.57 & 0.92 & 0.89 & 0.85 & 0.71 & 0.93 & 20.15 \\
Menger                                             & 0.49 & 0.01 & \underline{0.08} & 0.92 & 0.83 & \underline{0.46} & 0.63 & 0.41 & \underline{17.80} \\
Haantjes                                              & 0.49 & 0.67 & 0.49 & 0.92          & 0.66 & 22.43 &           &  & OOT \\
    \bottomrule
  \end{tabular}
\end{table}

In this ablation study, we observe that the choice of different types of Ricci curvature markedly influences both community detection quality and computational cost. For simplicity, {\textsf{PLRF}}s based on Ollivier's Ricci curvature, Lin-Lu-Yau's Ricci curvature, Forman's Ricci curvature,
Menger's Ricci curvature and Haantjes's Ricci curvature are shortened as Ollivier, Lin-Lu-Yau, Forman, Menger and Haantjes respectively.
The Ollivier consistently attains the highest NMI on the Karate (0.93) and Football (0.94) networks, as well as the top modularity on Football (0.92) and Facebook (0.95), demonstrating its ability to closely recover ground-truth partitions. However, this accuracy comes at the expense of runtime: Ollivier requires several seconds on small graphs and nearly half an hour on the Facebook network (1775.81s). In contrast, the Forman achieves a more favorable balance, delivering competitive NMI (0.49-0.92) and modularity (0.67-0.93) with run-times of subseconds, even on the largest graph (20.15s). The Lin-Lu-Yau offers the highest modularity on the Karate network (0.86) and remains under one second on Football, but its cost grows substantially on Facebook (1396.94s). Menger, although extremely fast, fails to produce a meaningful community structure (e.g., Q = 0.01 in Karate), and the Haantjes does not scale beyond medium-sized networks (OOT on Facebook).\\

In general, {\textsf{PLRF}} based on Ollivier's Ricci curvature is preferable when maximal detection quality is required and runtime constraints are less critical, while Forman's Ricci curvature-based {\textsf{PLRF}} stands out as the most practical choice for large-scale applications, striking an optimal balance between accuracy and efficiency. However, {\textsf{PLRF}} based on Lin-Lu-Yau's Ricci curvature may be selected when prioritizing modularity in small-graph analyses.

\section{Conclusion}\label{Sec 7}
The {\textsf PLRF} framework establishes a rigorous theoretical foundation for discrete Ricci
flows, ensuring global existence, uniqueness, and convergence while enabling effective
community detection. By bridging geometric analysis and graph theory, {\textsf PLRF} offers a powerful
tool for uncovering topological structures in complex systems. Its superior performance on
real-world datasets and theoretical guarantees position it as a valuable addition to the network
analysis toolkit. Moving forward, optimizing computational efficiency and exploring hybrid
geometric-learning approaches will further expand its utility in practical applications.


\end{document}